\documentclass[10pt]{article}
\usepackage[intlimits]{amsmath}
\usepackage{amssymb}
\usepackage{epsfig}
\usepackage{float}
\newtheorem{theorem}{Theorem}[section]

\newtheorem{definition}{Definition}[section]
\newtheorem{note}{Note}[section]
\newtheorem{proposition}{Proposition}[section]
\newtheorem{lemma}{Lemma}[section]

\newtheorem{rule-def}[theorem]{Rule}
\newtheorem{example}{Example}[section]

\begin{document}

\title{NONLINEAR THREE POINT SINGULAR BVPs : A CLASSIFICATION}
\author{Mandeep Singh$^a$\thanks{$^b$akverma@iitp.ac.in, amitkverma02@yahoo.co.in, $^a$ mandeep04may@yahoo.in}, Amit K. Verma$^b$
\\\small{\it{$^b$Department of Mathematics, IIT Patna, Bihta, Bihar--$801103$, India.}
}\\{\small\it{$^a$Department of Mathematics, DIT University, Dehradun, Uttarakhand-248009, India.}}}
\date{}
\maketitle
\begin{abstract}
We analyze the existence of unique solutions of the following class of nonlinear three point singular boundary value problems (SBVPs),
\begin{eqnarray*}\label{NL-Singular-P}
&&-(x^{\alpha} y'(x))'= x^{\alpha}f(x,y) ,\quad 0<x<1,\\
&&y'(0)=0,\quad y(1)=\delta y(\eta),
\end{eqnarray*}
where $\delta>0$, $0<\eta<1$ and $\alpha \geq 1$. This study shows some novel observations regarding the nature of the solution of the nonlinear three point SBVPs. We observe that when $sup\left(\partial f/\partial y\right)>0$ for $\alpha\in \cup_{n\in \mathbb{N}}\left(4n-1,4n+1\right)$ or $\alpha\in\{1,5,9,\cdots\}$ reverse ordered case occur. When $sup\left(\partial f/\partial y\right)>0$ for $\alpha\in \cup_{n\in \mathbb{N}}\left(4n-3,4n-1\right)$ or $\alpha\in\{3,7,11,\cdots\}$ and when $sup\left(\partial f/\partial y\right)<0$ for all $\alpha\geq 1$ well order case occur.
\end{abstract}


\section{\textbf{Introduction}}
Nonlinear singular boundary value problems, commonly arises due to physical symmetry. Specifically, if a physical law is expressed by a Dirichlet problem $\bigtriangleup y(P)=f(P,y(P)),$ and one is interested in planar, cylindrical or spherical geometries, then it is led to the differential equation
\begin{eqnarray}
&&\label{1}-y''-\frac{\alpha}{x}y'=f(x,y),\quad 0 < x < 1,
\end{eqnarray}
where $\alpha= 0,1,2$, respectively \cite{Russell-Shampine-SIAM-JNA-1975}. In modern science various real life problems have been converted into a mathematical model similar to (\ref{1}), where $\alpha\geq0$, and several researchers have studied existence and uniqueness of the solution of (\ref{1}), e.g., (see \cite{Russell-Shampine-SIAM-JNA-1975, P.L.Chamber-JCP-1952, SC-DOVER-1967, J.B.Keller-JRMA-1956,Chawla-shivkumar-JCAM-1987, AKV-NA-2011}).

If the condition at the end point $(x=1)$ depends on some interior point of the interval $(0,1)$, then the following three point nonlinear SBVP arises
\begin{eqnarray}
\label{SBVP-1}&&-(x^{\alpha} y'(x))'= x^{\alpha}f(x,y) ,\quad 0<x<1,\\
\label{SBVP-2}&&y'(0)=0,\quad y(1)=\delta y(\eta),
\end{eqnarray}
where $f:I\times R\rightarrow R,~I=[0,1],~0<\eta<1,~\delta>0$ and $\alpha\geq1$. Different analytical techniques are available for solvability of three point nonlinear BVPs, for $\alpha=0, 1$ and $2$ (see \cite{FL-MJ-XL-CL-GL-NA-2008, J.R.L.Webb-NA-RWA-2012, Gregus-Neumann-Arscott-JLMS-1971, R.Ma-JMAA-2001, JN-NATMA-1997, JH-BK-CCT-PAMC-2005, AKV-MS-JOMC-2015}. The results of this paper may lead to some new developments in this area.

Two point BVPs are studied extensively when upper and lower solutions are well ordered. But when upper and lower solutions are in the reversed order (see \cite{AKV-MS-JOMC-2015, Zhang-SJMA-1995, Cabada-BVP-2011, Coster-Habets-Elsevier-2006, MC-CDC-PH-AMC-2001,CA-PA-SL-AMC-2001} and the reference there in) lot of explorations are still left unattended for both two point and multi point boundary value problems. In this paper we have tried to analyze few of these gaps related to monotone iterative techniques.

If $f(x,y)$ is continuous and Lipschitz continuous in its domain, we derive the following monotone iterative scheme from the nonlinear singular three point BVP (\ref{SBVP-1})--(\ref{SBVP-2}),
\begin{eqnarray}\label{MIS}
-y_{n+1}''(x)-\frac{\alpha}{x}y_{n+1}'(x)-\lambda y_{n+1}(x)=f(x,y_n)-\lambda y_n(x),\quad y_{n+1}'(0)=0,\quad y_{n+1}(1)=\delta y_{n+1}(\eta),
\end{eqnarray}
and prove that solution exists and belongs to the class $C[0,1]\cap~C^2[0,1]$.

We also define upper and lower solution which are our initial guesses for the above iterative scheme.
\begin{definition}
If the functions $u_0,v_0 \in C^2[0,1]$ are defined as
\begin{eqnarray}
\label{US-1}&& -(x^{\alpha}u_0'(x))'\geq x^{\alpha}f(x,u_0),\quad 0<x<1;\\
\label{US-2}&&~~u_0'(0)=0,\quad u_0(1)\geq \delta u_0(\eta),\\
\label{LS-1}&& -(x^{\alpha}v_0'(x))'\leq x^{\alpha}f(x,v_0),\quad 0<x<1;\\
\label{LS-2}&&~~v_0'(0)=0,\quad v_0(1)\leq \delta v_0(\eta),
\end{eqnarray}
 then $u_0$ and $v_0$ are called upper and lower solutions of the nonlinear three point SBVPs (\ref{SBVP-1})--(\ref{SBVP-2}), respectively.
\end{definition}
The purpose of this paper is to prove existence of unique solution for the class of nonlinear three point SBVPs. We observe that depending on the values of $\alpha$ we arrive at well ordered and reversed order cases. This classification, we deduce does not exist in the literature to the best of our knowledge.

\begin{figure}[H]
\includegraphics[width=16cm,height=21cm]{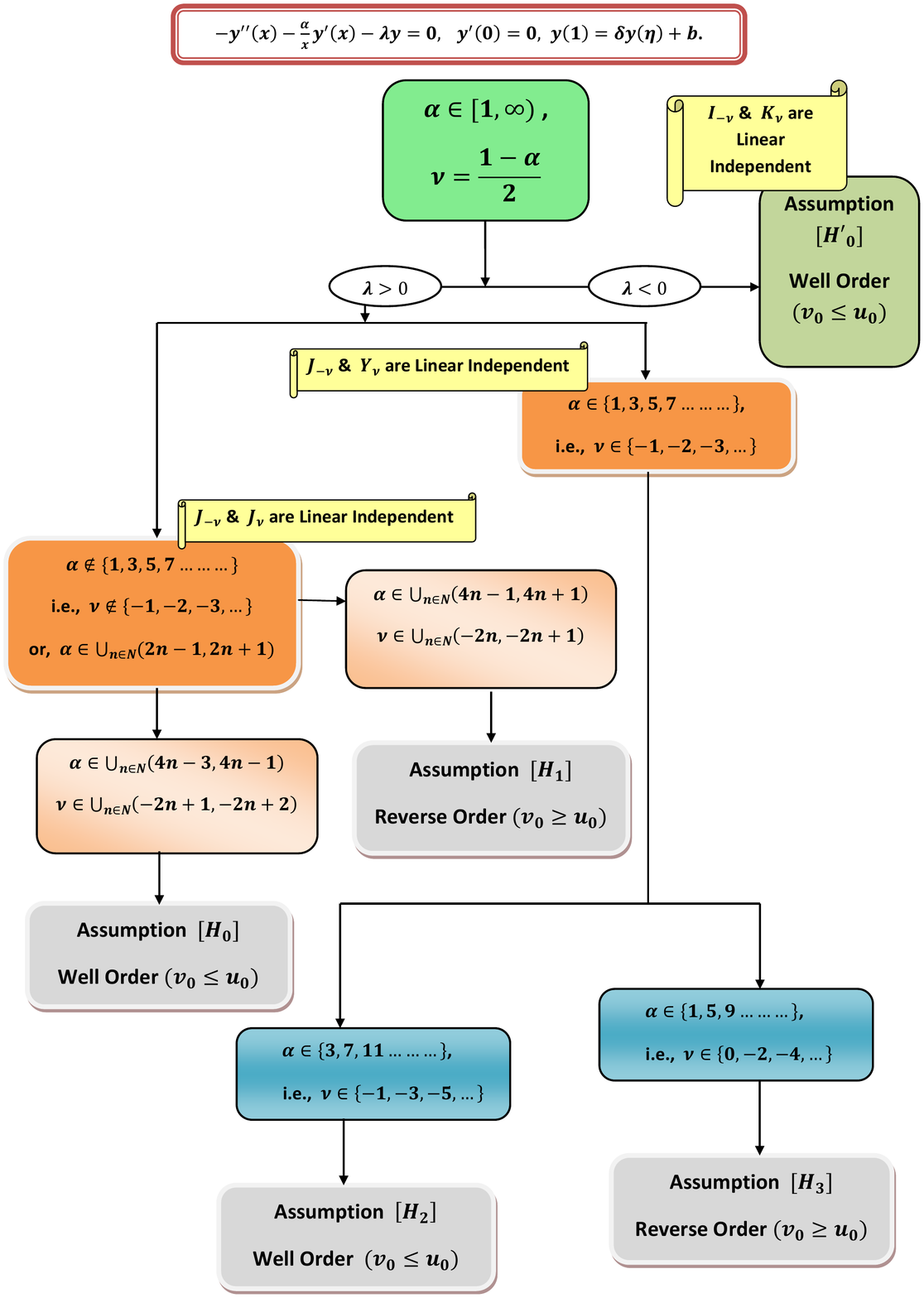}\\
\caption{Flow Chart}\label{Figure}
\end{figure}

This paper is organized in the following sections. Section \ref{TLC} we use Lommel's transformation to find out two linearly independent solutions in the terms of Bessel functions. Using these two linearly independent solutions Green's functions are constructed for different class of $\alpha$ (See Figure \ref{Figure}) in Section 3 and Section 4 states maximum  and anti-maximum  principles. Finally in Section 5 all these results are used to establish some new existence and uniqueness theorems. The sufficient conditions derived in this paper are verified for certain values which belongs to different classes of $\alpha$ in Section \ref{Example}.
\section{The Linear Case}\label{TLC}
The linear BVP corresponding to the nonlinear three point SBVPs (\ref{SBVP-1})--(\ref{SBVP-2}) is studied in this section. We consider the following inhomogeneous class of three point linear SBVPs,
\begin{eqnarray}
\label{Linear-SBVP}&&-(x^{\alpha} y'(x))'-\lambda x^{\alpha} y(x)=x^{\alpha} h(x),\quad 0<x<1,\\
&&\label{Linear-SBVP-B.C.-1}y'(0)=0,\quad y(1)=\delta y(\eta)+b,
\end{eqnarray}
where $h\in C(I)$ and $b$ is any constant. To solve the inhomogeneous system (\ref{Linear-SBVP})--(\ref{Linear-SBVP-B.C.-1}), we consider the corresponding homogeneous system
\begin{eqnarray}
\label{HP}&&-(x^{\alpha} y'(x))'-\lambda x^{\alpha} y(x)=0,\quad 0<x<1,\\
&&\label{HP-B.C.-1}y'(0)=0,\quad y(1)=\delta y(\eta).
\end{eqnarray}
Using Lommel's transformation (\S cf \cite{Chawla-shivkumar-JCAM-1987, A-Erdrlyi-McGraw-Hill-1953}) $z=\beta\zeta^{\gamma}$, $w= \zeta^{-a}v(\zeta)$, the standard Bessel's equation (\ref{2.4h}) is transformed into (\ref{2.5h})
\begin{eqnarray}
\label{2.4h}&&z^2\frac{d^2w}{dz^2}+z\frac{dw}{dz}+(z^2-\nu^2)w=0,\\
\label{2.5h}&&{\zeta}^2\frac{d^2v}{d{\zeta}^2}+{\zeta}(1-2a)\frac{dv}{d{\zeta}}+\left[(\beta \gamma \zeta^{\gamma})^2+(a^2-\nu^2\gamma^2)\right]v=0.
\end{eqnarray}
Now, by Lommel's Transformation the two linearly independent solutions of (\ref{2.5h}) are given by
\begin{eqnarray}
\label{2.6h}&&v_1(x)=\zeta^{a}w_1\left(\beta\zeta^{\gamma}\right),\quad v_2(x)= \zeta^{a}w_2\left(\beta\zeta^{\gamma}\right),
\end{eqnarray}
where $w_1(z)$ and $w_2(z)$ are two linearly independent solutions of Bessel's equation (\ref{2.4h}).\\
Now, if we set $\nu=a=\frac{1-\alpha}{2}$, $\gamma=1$, $\beta^2=\lambda$, then (\ref{2.5h}) reduces to (\ref{HP}) and hence we can obtained the two linearly independent solutions of (\ref{HP}) in terms of $w_1(z)$ and $w_2(z)$. A solution of (\ref{HP}) which is bounded in the neighborhood of the origin (except for a multiplicative constant) given by
$x^{\nu} J_{-\nu}\left(x \sqrt{\lambda}\right),~\mbox{if $\lambda>0$}~ \mathrm{and}~ x^{\nu} I_{-\nu}\left( x  \sqrt{|\lambda|}\right), ~\mbox{if $\lambda<0$}$.
\begin{note}
In this paper $J_{-\nu}$, $Y_{\nu}$ are Bessel functions of first and second kind and $I_{-\nu}$ and $K_{\nu}$ are Modified Bessel functions of first and second kind.
\end{note}
\section{Green's function} On the basis of sign of $\lambda$ and values of $\alpha$, we divide into the following cases.
\subsubsection{Case I: When $\lambda>0$ and $\alpha\notin \{1,3,5,\cdots\infty\}$.}
Suppose that
\begin{itemize}
\item [$(H_0):$]  $0<\lambda<j^{2}_{\nu,1}$,~~$0<\delta<1$,~~$\delta  \eta ^{\nu } J_{\nu }\left(\eta  \sqrt{\lambda }\right)-J_{\nu }\left(\sqrt{\lambda }\right)\geq 0$ and $J_{-\nu }\left(\sqrt{\lambda }\right)-\delta  \eta ^{\nu } J_{-\nu }\left(\eta  \sqrt{\lambda }\right)>0$,\\
      for $\alpha\in\bigcup_{n\in N} (4n-3,4n-1)$, and
\item [$(H_1):$]  $0<\lambda<j^{2}_{\nu,1}$,~~$\delta\geq1$,~~$\delta  \eta ^{\nu } J_{\nu }\left(\eta  \sqrt{\lambda }\right)-J_{\nu }\left(\sqrt{\lambda }\right)\leq 0$ and $J_{-\nu }\left(\sqrt{\lambda }\right)-\delta  \eta ^{\nu } J_{-\nu }\left(\eta  \sqrt{\lambda }\right)<0$,\\
     for $\alpha\in\bigcup_{n\in N} (4n-1,4n+1)$,
\end{itemize}
where $j^{2}_{\nu,1}$ is the first zero of $J_{\nu}(x)$. We can easily check the validation of $(H_0)$ and $(H_1)$ for $\alpha\in\bigcup_{n\in N} (4n-3,4n-1)$ and $\alpha\in\bigcup_{n\in N} (4n-1,4n+1)$ respectively.\\

Next two lemmas help us to define the sign of Green's function.

\pagebreak
\pagebreak

\begin{lemma}\label{LemmaJ-J-1} For $0<\lambda<\ j_{\nu,1}^2$, the Bessel functions of first kind $(J_{\nu}$ and $J_{-\nu})$ satisfy the following inequality
\begin{eqnarray*}
r^\nu \left(J_{-\nu}\left(s \sqrt{\lambda }\right)J_{\nu}\left(r \sqrt{\lambda }\right)-J_{\nu}\left(s \sqrt{\lambda }\right)J_{-\nu}\left(r \sqrt{\lambda }\right)\right)\geq 0,~~ 0< r\leq s\leq 1,
\end{eqnarray*}
where $\nu=\frac{1-\alpha}{2}$ and $\alpha\in\bigcup_{n\in N} (4n-3,4n-1)$.
\end{lemma}

\textbf{Proof.} Suppose
$$\widetilde{\Phi}(s,r)=r^\nu \left(J_{-\nu}\left(s \sqrt{\lambda }\right)J_{\nu}\left(r \sqrt{\lambda }\right)-J_{\nu}\left(s \sqrt{\lambda }\right)J_{-\nu}\left(r \sqrt{\lambda }\right)\right),$$ and let $s = s_0\in[0,1]$ be fixed. As
\begin{eqnarray*}
J_{-\nu}\left(s_0 \sqrt{\lambda }\right)J_{-1+\nu}\left(r \sqrt{\lambda }\right)+J_{\nu}\left(s_0 \sqrt{\lambda }\right)J_{1-\nu}\left(r \sqrt{\lambda }\right)&\leq&J_{-\nu}\left(s_0 \sqrt{\lambda }\right)J_{-1+\nu}\left(s_0 \sqrt{\lambda }\right)+J_{\nu}\left(s_0 \sqrt{\lambda }\right)J_{1-\nu}\left(x_0 \sqrt{\lambda }\right)\\
&=&\frac{2 \sin{\nu \pi}}{\pi s_0 \sqrt{\lambda}}\leq 0,
\end{eqnarray*}
for $r\leq s_0$. Now making use of the above inequalitiy, we can easily show that $\widetilde{\Phi}(s_0,t)$ is a non-increasing function of $r$. As $\widetilde{\Phi}(s_0,s_0)=0$ at $r=s_0$, which implies that $\widetilde{\Phi}(s_0,r)\geq 0,~~~\forall~~ r\leq s_0$. But as $s_0$ may have any value in $[0,1]$ therefore $\widetilde{\Phi}(s,r)\geq 0,~~\forall~~ 0< r\leq s\leq 1$.
\begin{lemma}\label{LemmaJ-J-2} For $0<\lambda<\ j_{\nu,1}^2$, the Bessel functions of first kind $(J_{\nu}$ and $J_{-\nu})$ satisfy the following inequality
\begin{eqnarray*}
r^\nu \left(J_{-\nu}\left(s \sqrt{\lambda }\right)J_{\nu}\left(r \sqrt{\lambda }\right)-J_{\nu}\left(s \sqrt{\lambda }\right)J_{-\nu}\left(r \sqrt{\lambda }\right)\right)\leq 0,~~ 0< r\leq s\leq 1,
\end{eqnarray*}
where $\nu=\frac{1-\alpha}{2}$ and $\alpha\in\bigcup_{n\in N} (4n-1,4n+1)$.
\end{lemma}

\textbf{Proof.} Proof follows the same analysis as we do in Lemma \ref{LemmaJ-J-1}, with the Bessel functions inequality
\begin{eqnarray*}
J_{-\nu}\left(s_0 \sqrt{\lambda }\right)J_{-1+\nu}\left(r \sqrt{\lambda }\right)+J_{\nu}\left(s_0 \sqrt{\lambda }\right)J_{1-\nu}\left(r \sqrt{\lambda }\right)&\geq&J_{-\nu}\left(s_0 \sqrt{\lambda }\right)J_{-1+\nu}\left(s_0 \sqrt{\lambda }\right)+J_{\nu}\left(s_0 \sqrt{\lambda }\right)J_{1-\nu}\left(s_0 \sqrt{\lambda }\right)\\
&=&\frac{2 \sin{\nu \pi}}{\pi s_0 \sqrt{\lambda}}\geq 0,
\end{eqnarray*}
for $r \leq s_0$.
\begin{lemma}\label{Lemma2.1}
For the linear three point SBVPs (\ref{HP})--(\ref{HP-B.C.-1}), where $\alpha\notin \{1,3,5,\cdots\infty\}$, the Green's function is given by
\begin{eqnarray*}
G(x,t)=\left\{
  \begin{array}{ll}
   \frac{\pi  \csc (\pi  \nu ) t^{\nu } x^{\nu } J_{-\nu }\left(x \sqrt{\lambda }\right) \left(\delta  \eta ^{\nu } \left(J_{\nu }\left(\eta  \sqrt{\lambda }\right) J_{-\nu }\left(t \sqrt{\lambda }\right)-J_{-\nu }\left(\eta  \sqrt{\lambda }\right) J_{\nu }\left(t \sqrt{\lambda }\right)\right)+\left(J_{-\nu }\left(\sqrt{\lambda }\right) J_{\nu }\left(t \sqrt{\lambda }\right)-J_{\nu }\left(\sqrt{\lambda }\right) J_{-\nu }\left(t \sqrt{\lambda }\right)\right)\right)}{2 \left(J_{-\nu }\left(\sqrt{\lambda }\right)- \delta  \eta ^{\nu } J_{-\nu }\left(\eta  \sqrt{\lambda }\right)\right)}, & 0 \leq x\leq t\leq \eta; \\
  \frac{1}{2} \pi  \csc (\pi  \nu ) t^{\nu } x^{\nu } J_{-\nu }\left(t \sqrt{\lambda }\right) \left(\frac{J_{-\nu }\left(x \sqrt{\lambda }\right) \left(\delta  \eta ^{\nu } J_{\nu }\left(\eta  \sqrt{\lambda }\right)-J_{\nu }\left(\sqrt{\lambda }\right)\right)}{J_{-\nu }\left(\sqrt{\lambda }\right)-\delta  \eta ^{\nu } J_{-\nu }\left(\eta  \sqrt{\lambda }\right)}+J_{\nu }\left(x \sqrt{\lambda }\right)\right), & t \leq x,~ t \leq \eta;\\
  \frac{\pi  \csc (\pi  \nu ) t^{\nu } x^{\nu } \left(J_{-\nu }\left(\sqrt{\lambda }\right) J_{\nu }\left(t \sqrt{\lambda }\right)-J_{\nu }\left(\sqrt{\lambda }\right) J_{-\nu }\left(t \sqrt{\lambda }\right)\right) J_{-\nu }\left(x \sqrt{\lambda }\right)}{2 \left(J_{-\nu }\left(\sqrt{\lambda }\right)- \delta  \eta ^{\nu } J_{-\nu }\left(\eta  \sqrt{\lambda }\right)\right)}, & x\leq t, ~\eta \leq t;\\
   \frac{1}{2} \pi  \csc (\pi  \nu ) t^{\nu } x^{\nu } \left(\frac{J_{-\nu }\left(x \sqrt{\lambda }\right) \left(\delta  \eta ^{\nu } J_{-\nu }\left(\eta  \sqrt{\lambda }\right) J_{\nu }\left(t \sqrt{\lambda }\right)-J_{\nu }\left(\sqrt{\lambda }\right) J_{-\nu }\left(t \sqrt{\lambda }\right)\right)}{J_{-\nu }\left(\sqrt{\lambda }\right)-\delta  \eta ^{\nu } J_{-\nu }\left(\eta  \sqrt{\lambda }\right)}+J_{-\nu }\left(t \sqrt{\lambda }\right) J_{\nu }\left(x \sqrt{\lambda }\right)\right), & \eta \leq t\leq x\leq 1.
  \end{array}
\right.
\end{eqnarray*}
If $H_0$ (or $H_1$) holds then $G(x,t)\leq0$ (or  $G(x,t)\geq0$).
\end{lemma}

\textbf{Proof.} We define Green's function as
\begin{eqnarray*}
G(x,t)=\left\{
  \begin{array}{ll}
   a_1~ x^{\nu } J_{-\nu }\left(x \sqrt{\lambda }\right)+a_2~ x^{\nu } J_{\nu }\left(x \sqrt{\lambda }\right) , & {0 \leq x\leq t\leq \eta;} \\
   a_3~ x^{\nu } J_{-\nu }\left(x \sqrt{\lambda }\right)+a_4~ x^{\nu } J_{\nu }\left(x \sqrt{\lambda }\right), & {t \leq x,~ t \leq \eta;} \\
   a_5~ x^{\nu } J_{-\nu }\left(x \sqrt{\lambda }\right)+a_6~ x^{\nu } J_{\nu }\left(x \sqrt{\lambda }\right), & {x \leq t,~ \eta \leq t;} \\
   a_7~ x^{\nu } J_{-\nu }\left(x \sqrt{\lambda }\right)+a_8~ x^{\nu } J_{\nu }\left(x \sqrt{\lambda }\right), & {\eta \leq t\leq x\leq 1}.
  \end{array}
\right.
\end{eqnarray*}
Using the following properties of the Green's function, for any $t \in [0,\eta]$ and boundary conditions, we have
\begin{eqnarray*}
&& a_1 t^{\nu}J_{-\nu}\left(t \sqrt{\lambda }\right)+ a_2 t^{\nu} J_{\nu}\left(t\sqrt{\lambda}\right)=a_3  t^{\nu}J_{-\nu}\left(t \sqrt{\lambda }\right)+ a_4 t^{\nu} J_{\nu}\left(t\sqrt{\lambda}\right),\\
&&-a_1 t^{\nu}\sqrt{\lambda }J_{1-\nu}\left(t \sqrt{\lambda }\right)+a_2 t^{\nu} \sqrt{\lambda} J_{\nu-1}\left(t\sqrt{\lambda}\right)+a_3 t^{\nu}\sqrt{\lambda }J_{1-\nu}\left(t \sqrt{\lambda }\right)-a_4 t^{\nu} \sqrt{\lambda} J_{\nu-1}\left(t\sqrt{\lambda}\right) = -\frac{1}{t^{1-2\nu}}
\end{eqnarray*}
and the following system of equations is derived
$$
\left(
  \begin{array}{cccc}
1 & 0 & -1 & 0\\
0 & 1 & 0 & -1\\
0 & 1 & 0 & 0\\
0 & 0 & J_{-\nu}\left( \sqrt{\lambda }\right)-\delta \eta^{\nu}J_{-\nu}\left(\eta \sqrt{\lambda }\right) & J_{\nu}\left(\sqrt{\lambda}\right)-\delta\eta^{\nu}J_{\nu}\left(\eta \sqrt{\lambda }\right)\\
  \end{array}
\right)
\left(
  \begin{array}{c}
    a_1 \\
    a_2 \\
    a_3 \\
    a_4 \\
  \end{array}
\right)=
\left(
  \begin{array}{c}
    \frac{\pi t^{\nu } J_{\nu }\left(t \sqrt{\lambda }\right)}{2 \sin{\nu \pi}} \\
    -\frac{\pi t^{\nu } J_{-\nu }\left(t \sqrt{\lambda }\right)}{2\sin{\nu\pi}} \\
      0 \\
      0 \\
 \end{array}
\right).$$
Solution of above system gives,
\begin{eqnarray*}
&&a_1=\frac{\pi  \csc (\pi  \nu ) t^{\nu } \left(J_{-\nu }\left(t \sqrt{\lambda }\right) \left(\delta  \eta ^{\nu } J_{\nu }\left(\eta  \sqrt{\lambda }\right)-J_{\nu }\left(\sqrt{\lambda }\right)\right)+J_{\nu }\left(t \sqrt{\lambda }\right) \left(J_{-\nu }\left(\sqrt{\lambda }\right)-\delta  \eta ^{\nu } J_{-\nu }\left(\eta  \sqrt{\lambda }\right)\right)\right)}{2 \left(J_{-\nu }\left(\sqrt{\lambda }\right)- \delta  \eta ^{\nu } J_{-\nu }\left(\eta  \sqrt{\lambda }\right)\right)},\\
&&a_2=0,\\
&&a_3=\frac{\pi  \csc (\pi  \nu ) t^{\nu } J_{-\nu }\left(t \sqrt{\lambda }\right) \left(\delta  \eta ^{\nu } J_{\nu }\left(\eta  \sqrt{\lambda }\right)-J_{\nu }\left(\sqrt{\lambda }\right)\right)}{2 \left(J_{-\nu }\left(\sqrt{\lambda }\right)- \delta  \eta ^{\nu } J_{-\nu }\left(\eta  \sqrt{\lambda }\right)\right)},\\
&&a_4=\frac{1}{2} \pi  \csc (\pi  \nu ) t^{\nu } J_{-\nu }\left(t \sqrt{\lambda }\right).
\end{eqnarray*}
Similarly for any $t \in [\eta,1]$, we have
\begin{eqnarray*}
&& a_5 t^{\nu}J_{-\nu}\left(t \sqrt{\lambda }\right)+ a_6 t^{\nu} J_{\nu}\left(t\sqrt{\lambda}\right)=a_7  t^{\nu}J_{-\nu}\left(t \sqrt{\lambda }\right)+ a_8 t^{\nu} J_{\nu}\left(t\sqrt{\lambda}\right),\\
&&-a_5 t^{\nu}\sqrt{\lambda }J_{1-\nu}\left(t \sqrt{\lambda }\right)+a_6 t^{\nu} \sqrt{\lambda} J_{\nu-1}\left(t\sqrt{\lambda}\right)+a_7 t^{\nu}\sqrt{\lambda }J_{1-\nu}\left(t \sqrt{\lambda }\right)-a_8 t^{\nu} \sqrt{\lambda} J_{\nu-1}\left(t\sqrt{\lambda}\right) = -\frac{1}{t^{1-2\nu}}.
\end{eqnarray*}
The above two equations and boundary conditions in $[\eta,1]$ gives
$$
\left(
  \begin{array}{cccc}
1 & 0 & -1 & 0\\
0 & 1 & 0 & -1\\
0 & 1 & 0 & 0\\
\delta \eta^{\nu}J_{-\nu}\left(\eta \sqrt{\lambda }\right) & \delta \eta^{\nu}J_{\nu}\left(\eta \sqrt{\lambda }\right) & -J_{-\nu}\left(\sqrt{\lambda }\right)  &  -J_{\nu}\left(\sqrt{\lambda }\right)\\
  \end{array}
\right)
\left(
  \begin{array}{c}
    a_5 \\
    a_6 \\
    a_7 \\
    a_8 \\
  \end{array}
\right)=
\left(
  \begin{array}{c}
    \frac{\pi t^{\nu } J_{\nu }\left(t \sqrt{\lambda }\right)}{2 \sin{\nu \pi}} \\
    -\frac{\pi t^{\nu } J_{-\nu }\left(t \sqrt{\lambda }\right)}{2\sin{\nu\pi}} \\
      0 \\
      0 \\
 \end{array}
\right).$$
By above four equation we have
\begin{eqnarray*}
&&a_5=\frac{\pi  \csc (\pi  \nu ) t^{\nu } \left(J_{-\nu }\left(\sqrt{\lambda }\right) J_{\nu }\left(t \sqrt{\lambda }\right)-J_{\nu }\left(\sqrt{\lambda }\right) J_{-\nu }\left(t \sqrt{\lambda }\right)\right)}{2 \left(J_{-\nu }\left(\sqrt{\lambda }\right)- \delta  \eta ^{\nu } J_{-\nu }\left(\eta  \sqrt{\lambda }\right)\right)},\\
&&a_6=0,\\
&&a_7=\frac{\pi  \csc (\pi  \nu ) t^{\nu } \left(\delta  \eta ^{\nu } J_{-\nu }\left(\eta  \sqrt{\lambda }\right) J_{\nu }\left(t \sqrt{\lambda }\right)-J_{\nu }\left(\sqrt{\lambda }\right) J_{-\nu }\left(t \sqrt{\lambda }\right)\right)}{2 \left(J_{-\nu }\left(\sqrt{\lambda }\right)- \delta  \eta ^{\nu } J_{-\nu }\left(\eta  \sqrt{\lambda }\right)\right)},\\
&&a_8=\frac{1}{2} \pi  \csc (\pi  \nu ) t^{\nu } J_{-\nu }\left(t \sqrt{\lambda }\right).
\end{eqnarray*}
This completes the construction of Green's function. Using $(H_0)$~(or $H_1$) and Lemma \ref{LemmaJ-J-1} (or Lemma \ref{LemmaJ-J-2}) we can easily verify that $G(x,t)\leq0$ (or  $G(x,t)\geq0$).
\subsubsection{Case II: When $\lambda>0$ and $\alpha\in \{1,3,5,\cdots\infty\}$.}
Suppose that
\begin{itemize}
\item [$(H_2):$]  $0<\lambda<y^{2}_{\nu,1}$,~~$0<\delta<1$,~~$\delta  \eta ^{\nu } Y_{\nu }\left(\eta  \sqrt{\lambda }\right)-Y_{\nu }\left(\sqrt{\lambda }\right)\geq 0$ and $J_{-\nu }\left(\sqrt{\lambda }\right)-\delta  \eta ^{\nu } J_{-\nu }\left(\eta  \sqrt{\lambda }\right)>0$, for~~$\alpha\in \{3,7,11,\cdots\}$, and
\item [$(H_3):$]  $0<\lambda<y^{2}_{\nu,1}$,~~$\delta\geq1$,~~$\delta  \eta ^{\nu } Y_{\nu }\left(\eta \sqrt{\lambda }\right)-Y_{\nu }\left(\sqrt{\lambda }\right)\leq 0$ and $J_{-\nu }\left(\sqrt{\lambda }\right)-\delta  \eta ^{\nu } J_{-\nu }\left(\eta \sqrt{\lambda }\right)<0$, for $\alpha\in\{1,5,9,\cdots\}$,

    where $y^{2}_{\nu,1}$ is the first zero of $Y_{\nu}(x)$.  It is easy to see that $(H_2)$ and $(H_3)$ can be satisfied for $\alpha\in \{3,7,11,\cdots\}$ and $\alpha\in \{1,5,9,\cdots\}$, respectively.
\end{itemize}

\begin{lemma}\label{LemmaJ-Y-1} For $0<\lambda<\ y_{\nu,1}^2$, the Bessel functions of first and second kind $(J_{-\nu}$ and $Y_{\nu})$ satisfy the following inequality
\begin{eqnarray*}
r^\nu\left(J_{-\nu}\left(s \sqrt{\lambda }\right)Y_{\nu}\left(r \sqrt{\lambda }\right)-Y_{\nu}\left(s \sqrt{\lambda }\right)J_{-\nu}\left(r \sqrt{\lambda }\right)\right)\geq 0,~~ 0 <r\leq s\leq 1,
\end{eqnarray*}
where $\nu=\frac{1-\alpha}{2}$ and $\alpha\in \{3,7,11,\cdots\}$.
\end{lemma}
\textbf{Proof.} Suppose
$$\widetilde{F}(s,r)=r^\nu\left(J_{-\nu}\left(s \sqrt{\lambda }\right)Y_{\nu}\left(r \sqrt{\lambda }\right)-Y_{\nu}\left(s \sqrt{\lambda }\right)J_{-\nu}\left(r \sqrt{\lambda }\right)\right),$$ and let $s = s_0\in[0,1]$ be fixed. Now as $J_{-\nu}\left(s_0 \sqrt{\lambda }\right)\geq J_{1-\nu}\left(r \sqrt{\lambda }\right)$ and $-Y_{\nu}\left(s_0 \sqrt{\lambda }\right) \geq Y_{-1+\nu}\left(r \sqrt{\lambda }\right)$ for $r \leq s_0$, when $0<\lambda\leq\ y_{v,1}^2$, where $\alpha\in \{3,7,11,\cdots\}$. So with the help of these inequalities we can easily show that $\widetilde{F}(s_0,r)$ is an non-increasing function of $r$. As $\widetilde{F}(s_0,s_0)=0$ at $r=s_0$, which implies that $\widetilde{F}(s_0,r)\geq 0,~~~\forall~~ r\leq s_0$. But as $s_0$ may take any value in $[0,1]$ therefore $\widetilde{F}(s,r)\geq 0,~~\forall~~ 0<r\leq s\leq 1$.

\begin{lemma}\label{LemmaJ-Y-2} For $0<\lambda<\ y_{\nu,1}^2$, the Bessel functions of first and second kind $(J_{-\nu}$ and $Y_{\nu})$ satisfy the following inequality
\begin{eqnarray*}
r^\nu\left(J_{-\nu}\left(s \sqrt{\lambda }\right)Y_{\nu}\left(r \sqrt{\lambda }\right)-Y_{\nu}\left(s \sqrt{\lambda }\right)J_{-\nu}\left(r \sqrt{\lambda }\right)\right)\leq 0,~~ 0< r\leq s\leq 1,
\end{eqnarray*}
where $\nu=\frac{1-\alpha}{2}$ and $\alpha\in\{1,5,9,\cdots\}$.
\end{lemma}

\textbf{Proof:} By using the inequalities $J_{-\nu}\left(s_0 \sqrt{\lambda }\right)\geq J_{1-\nu}\left(r \sqrt{\lambda }\right)$ and $Y_{\nu}\left(s_0 \sqrt{\lambda }\right) \geq -Y_{-1+\nu}\left(r \sqrt{\lambda }\right)$ for $r \leq s_0$, when $0<\lambda\leq\ y_{v,1}^2$, where $\alpha\in\{1,5,9,\cdots\}$, we can prove this lemma as we did in Lemma \ref{LemmaJ-Y-1}.

\begin{lemma}\label{Lemma2.2}
For the linear three point SBVPs (\ref{HP})--(\ref{HP-B.C.-1}), where $\alpha\in \{1,3,5,\cdots\infty\}$, the Green's function is given by
\begin{eqnarray*}
G(x,t)=\left\{
  \begin{array}{ll}
   \frac{\pi  \sec (\pi  \nu ) t^{\nu } x^{\nu } J_{-\nu }\left(x \sqrt{\lambda }\right) \left(J_{-\nu }\left(t \sqrt{\lambda }\right) \left(\delta  \eta ^{\nu } Y_{\nu }\left(\eta  \sqrt{\lambda }\right)-Y_{\nu }\left(\sqrt{\lambda }\right)\right)+Y_{\nu }\left(t \sqrt{\lambda }\right) \left(J_{-\nu }\left(\sqrt{\lambda }\right)-\delta  \eta ^{\nu } J_{-\nu }\left(\eta  \sqrt{\lambda }\right)\right)\right)}{2\left( J_{-\nu }\left(\sqrt{\lambda }\right)- \delta  \eta ^{\nu } J_{-\nu }\left(\eta  \sqrt{\lambda }\right)\right)}, & 0 \leq x\leq t\leq \eta; \\
 \frac{1}{2} \pi  \sec (\pi  \nu ) t^{\nu } x^{\nu } J_{-\nu }\left(t \sqrt{\lambda }\right) \left(\frac{J_{-\nu }\left(x \sqrt{\lambda }\right) \left(\delta  \eta ^{\nu } Y_{\nu }\left(\eta  \sqrt{\lambda }\right)-Y_{\nu }\left(\sqrt{\lambda }\right)\right)}{J_{-\nu }\left(\sqrt{\lambda }\right)-\delta  \eta ^{\nu } J_{-\nu }\left(\eta  \sqrt{\lambda }\right)}+Y_{\nu }\left(x \sqrt{\lambda }\right)\right), & t \leq x,~ t \leq \eta;\\
  \frac{\pi  \sec (\pi  \nu ) t^{\nu } x^{\nu } J_{-\nu }\left(x \sqrt{\lambda }\right) \left(J_{-\nu }\left(\sqrt{\lambda }\right) Y_{\nu }\left(t \sqrt{\lambda }\right)-Y_{\nu }\left(\sqrt{\lambda }\right) J_{-\nu }\left(t \sqrt{\lambda }\right)\right)}{2\left( J_{-\nu }\left(\sqrt{\lambda }\right)- \delta  \eta ^{\nu } J_{-\nu }\left(\eta  \sqrt{\lambda }\right)\right)}, & x\leq t, ~\eta \leq t;\\
  \frac{1}{2} \pi  \sec (\pi  \nu ) t^{\nu } x^{\nu } \left(\frac{J_{-\nu }\left(x \sqrt{\lambda }\right) \left(\delta  \eta ^{\nu } J_{-\nu }\left(\eta  \sqrt{\lambda }\right) Y_{\nu }\left(t \sqrt{\lambda }\right)-Y_{\nu }\left(\sqrt{\lambda }\right) J_{-\nu }\left(t \sqrt{\lambda }\right)\right)}{J_{-\nu }\left(\sqrt{\lambda }\right)-\delta  \eta ^{\nu } J_{-\nu }\left(\eta  \sqrt{\lambda }\right)}+J_{-\nu }\left(t \sqrt{\lambda }\right) Y_{\nu }\left(x \sqrt{\lambda }\right)\right), & \eta \leq t\leq x\leq 1.
  \end{array}
\right.
\end{eqnarray*}
If $H_2$ (or $H_3$) holds then $G(x,t)\leq0$ (or  $G(x,t)\geq0$).
\end{lemma}
\textbf{Proof:} The construction of Green's function follows the same analysis as we do in Lemma \ref{Lemma2.1}. Using the assumption $(H_2)$~(or $H_3$) and Lemma \ref{LemmaJ-Y-1} (or Lemma \ref{LemmaJ-Y-2}) we can easily verify that $G(x,t)\leq0$ (or  $G(x,t)\geq0$).

Now we state Lemmas \ref{Lemma2.3}, \ref{Lemma2.5} and \ref{Lemma2.6} and we omit proof.
\begin{lemma}\label{Lemma2.3}
If $y\in C^2(I)$ is a solution of inhomogeneous linear three point SBVPs (\ref{Linear-SBVP})--(\ref{Linear-SBVP-B.C.-1}), then
\begin{eqnarray}
\label{2.5}&&y(x)= \frac{b ~ x^{\nu}J_{-\nu}\left(x \sqrt{\lambda }\right)}{J_{-\nu }\left(\sqrt{\lambda }\right)-\delta  \eta ^{\nu } J_{-\nu }\left(\eta  \sqrt{\lambda }\right)}- {\int_0}^1{t^{\alpha}~ G(x,t) h(t)dt}.
\end{eqnarray}
\end{lemma}
\subsubsection{Case III: When $\lambda<0$.}
Suppose that
\begin{itemize}
\item [$(H'_0):$] $\delta>0$,~~$K_{\nu }\left(\sqrt{|\lambda| }\right)-\delta  \eta ^{\nu } K_{\nu }\left(\eta  \sqrt{|\lambda|  }\right)\leq 0$ and $I_{-\nu }\left(\sqrt{|\lambda|  }\right)-\delta  \eta ^{\nu } I_{-\nu }\left(\eta  \sqrt{|\lambda|  }\right)>0$ for~~$\alpha\in[1,\infty)$.
\end{itemize}

\begin{lemma}\label{LemmaI-K}
For $\lambda< 0$, the modified Bessel functions of first and second kind $(I_{-\nu}$ and $K_{\nu})$ satisfy the following inequality
\begin{eqnarray*}
r^\nu\left(K_{\nu }\left(s \sqrt{ |\lambda|  }\right) I_{-\nu }\left(r \sqrt{|\lambda|  }\right)-I_{-\nu }\left(s \sqrt{|\lambda|  }\right) K_{\nu }\left(r \sqrt{|\lambda|  }\right)\right)\leq 0,~~~ 0 <r\leq s\leq 1
\end{eqnarray*}
where $\nu=\frac{1-\alpha}{2}$ and $\alpha\in[1,\infty)$.
\end{lemma}
\textbf{Proof:}
Suppose
$$\widetilde{F_1}(s,r)=r^\nu\left(K_{\nu }\left(s \sqrt{ |\lambda|  }\right) I_{-\nu }\left(r \sqrt{|\lambda|  }\right)-I_{-\nu }\left(s \sqrt{|\lambda|  }\right) K_{\nu }\left(r \sqrt{|\lambda|  }\right)\right),$$ and further assume that $s = s_0\in[0,1]$ be fixed. Now we can easily show that the function $\widetilde{F_1}(s_0,r)$ will be non-decreasing for $r$ for all $\alpha\in[1,\infty)$. At $r=s_0$, $\widetilde{F_1}(s_0,s_0)=0$, i.e., $\widetilde{F_1}(s_0,r)\leq 0,~~\forall~~r\leq s_0$. But as $s_0$ may have any value in $[0,1]$ therefore $\widetilde{F}(s,r)\leq 0,~\forall~~ 0<r\leq s\leq 1$.
\begin{lemma}\label{Lemma2.5}
For the following linear three point SBVPs (\ref{HP})--(\ref{HP-B.C.-1}), where $\alpha\in [1,\infty)$, the Green's function is given by
\begin{eqnarray*}
G(x,t)=\left\{
  \begin{array}{ll}
   \frac{t^{\nu } x^{\nu } I_{-\nu }\left(x \sqrt{|\lambda|  }\right) \left(I_{-\nu }\left(t \sqrt{|\lambda|  }\right) \left(K_{\nu }\left(\sqrt{|\lambda|  }\right)-\delta  \eta ^{\nu } K_{\nu }\left(\eta  \sqrt{|\lambda|  }\right)\right)-K_{\nu }\left(t \sqrt{|\lambda|  }\right) \left(I_{-\nu }\left(\sqrt{|\lambda|  }\right)-\delta  \eta ^{\nu } I_{-\nu }\left(\eta  \sqrt{|\lambda|  }\right)\right)\right)}{I_{-\nu }\left(\sqrt{|\lambda|  }\right)-\delta  \eta ^{\nu } I_{-\nu }\left(\eta  \sqrt{|\lambda|  }\right)},~~~0 \leq x\leq t\leq \eta; \\
  \frac{t^{\nu } x^{\nu } I_{-\nu }\left(t \sqrt{|\lambda|  }\right) \left(I_{-\nu }\left(x \sqrt{|\lambda|  }\right) \left(K_{\nu }\left(\sqrt{|\lambda|  }\right)-\delta  \eta ^{\nu } K_{\nu }\left(\eta  \sqrt{|\lambda|  }\right)\right)-K_{\nu }\left(x \sqrt{|\lambda|  }\right) \left(I_{-\nu }\left(\sqrt{|\lambda|  }\right)-\delta  \eta ^{\nu } I_{-\nu }\left(\eta  \sqrt{|\lambda|  }\right)\right)\right)}{I_{-\nu }\left(\sqrt{|\lambda|  }\right)-\delta  \eta ^{\nu } I_{-\nu }\left(\eta  \sqrt{|\lambda|  }\right)},~~~ t \leq x,~ t \leq \eta;\\
  \frac{t^{\nu } x^{\nu } I_{-\nu }\left(x \sqrt{|\lambda|  }\right) \left(K_{\nu }\left(\sqrt{|\lambda|  }\right) I_{-\nu }\left(t \sqrt{|\lambda|  }\right)-I_{-\nu }\left(\sqrt{|\lambda|  }\right) K_{\nu }\left(t \sqrt{|\lambda|  }\right)\right)}{I_{-\nu }\left(\sqrt{|\lambda|  }\right)-\delta  \eta ^{\nu } I_{-\nu }\left(\eta  \sqrt{|\lambda|  }\right)}, ~~~~~~~~~~~~~~~~~~~~~~~~~~~~~~~~~~~~~~~~~~~~ x\leq t, ~\eta \leq t;\\
    \frac{t^{\nu } x^{\nu } \left(I_{-\nu }\left(x \sqrt{|\lambda|  }\right) \left(K_{\nu }\left(\sqrt{|\lambda|  }\right) I_{-\nu }\left(t \sqrt{|\lambda|  }\right)-\delta  \eta ^{\nu } I_{-\nu }\left(\eta  \sqrt{|\lambda|  }\right) K_{\nu }\left(t \sqrt{|\lambda|  }\right)\right)-I_{-\nu }\left(t \sqrt{|\lambda|  }\right) K_{\nu }\left(x \sqrt{|\lambda|  }\right) \left(I_{-\nu }\left(\sqrt{|\lambda|  }\right)-\delta  \eta ^{\nu } I_{-\nu }\left(\eta  \sqrt{|\lambda|  }\right)\right)\right)}{I_{-\nu }\left(\sqrt{|\lambda|  }\right)-\delta  \eta ^{\nu } I_{-\nu }\left(\eta  \sqrt{|\lambda|  }\right)}, \\~~~~~~~~~~~~~~~~~~~~~~~~~~~~~~~~~~~~~~~~~~~~~~~~~~~~~~~~~~~~~~~~~~~~~~~~~~~~~~~~~~~~~~~~~~~~~~~~~~~~~~~~~~~~~~~~~ \eta \leq t\leq x\leq 1,
  \end{array}
\right.
\end{eqnarray*}
and if $(H'_0)$ holds then $G(x,t)\leq0$.
\end{lemma}
\begin{lemma}\label{Lemma2.6}
If $y\in C^2(I)$ is a solution of inhomogeneous linear three point SBVPs (\ref{Linear-SBVP})-(\ref{Linear-SBVP-B.C.-1}) then
\begin{eqnarray*}
&&y(x)= \frac{b ~ x^{\nu}I_{-\nu}\left(x \sqrt{|\lambda|  }\right)}{I_{-\nu }\left(\sqrt{|\lambda|  }\right)-\delta  \eta ^{\nu } I_{-\nu }\left(\eta  \sqrt{|\lambda|  }\right)}- {\int_0}^1{t^{\alpha}~ G(x,t) h(t)dt}.
\end{eqnarray*}
\end{lemma}

\section{Maximum and anti-maximum principles for linear three point SBVPs}
The constant sign of Green's function results into the following anti-Maximum and maximum principles.
\begin{proposition}\label{Proposition2.1}
\textbf{$($Anti-maximum principle$)$}\\
Assume $\lambda>0$ and $(H_1)$ or $(H_3)$ holds, and $y\in C^2(I)$ satisfies
\begin{eqnarray*}
&&-(x^{\alpha} y'(x))'-\lambda x^{\alpha} y(x)\geq0,\quad 0<x<1,\\
&&y'(0)=0,\quad y(1)\geq \delta y(\eta).
\end{eqnarray*}
Then $y(x)\leq0$, $\forall x\in [0,1]$.
\end{proposition}
\begin{proposition}\label{Proposition2.2}
\textbf{$($Maximum principle$)$}
\begin{itemize}
\item[$(Max_1)$] Assume $\lambda>0$ and $(H_0)$ or $(H_2)$ holds, and $y\in C^2(I)$ satisfies
\begin{eqnarray*}
&&-(x^{\alpha} y'(x))'-\lambda x^{\alpha} y(x)\geq0,\quad 0<x<1,\\
&&y'(0)=0,\quad y(1)\geq \delta y(\eta).
\end{eqnarray*}
Then $y(x)\geq0$, $\forall x\in [0,1]$.
\item[$(Max_2)$]  Assume $\lambda<0$, $(H'_0)$ holds and $y\in C^2(I)$ satisfies
\begin{eqnarray*}
&&-(x^{\alpha} y'(x))'-\lambda x^{\alpha} y(x)\geq0,\quad 0<x<1,\\
&&y'(0)=0,\quad y(1)\geq \delta y(\eta).
\end{eqnarray*}
Then $y(x)\geq0$, $\forall x\in [0,1]$.
\end{itemize}
\end{proposition}
\section{\textbf{Existence results for nonlinear three point SBVP}}
The study of existence results for nonlinear three point SBVP is discussed in this section. On the basis of anti-maximum and maximum principles, we divide it into  the following two subsections.
\subsubsection{Reverse ordered upper $(u_0)$ and lower $(v_0)$ solutions $(u_0\leq v_0).$}
\begin{theorem}\label{Theorem3.1}
Assume that
\begin{itemize}
\item[$(R_1)$] there exist $u_0$, $v_0$ in $C^2[0,1]$ such that $u_0\leq v_0$, where $u_0$ satisfies (\ref{US-1})--(\ref{US-2}) and $v_0$ satisfies (\ref{LS-1})--(\ref{LS-2});
\item[$(R_2)$] the function $f: D \rightarrow R$ is continuous on $D :=\{(x,y) \in [0,1] \times R :u_0 \leq y \leq v_0\}$;
\item[$(R_3)$] there exist $M\geq 0$ such that for all $(x,y),(x,w)\in D$
\begin{eqnarray*}
y\leq w \Longrightarrow  f(x,w)- f(x,y) \leq M_1(w - y),
\end{eqnarray*}
\item[$(R_4)$] there exist a constant $\lambda>0$ such that $M_1-\lambda\leq0$ and $H_1$ or $H_3$ holds;
\end{itemize}
then the nonlinear three point SBVP (\ref{SBVP-1})--(\ref{SBVP-2}) has at least one solution in the region $D$. Sequence $\{u_n\}$ generated by equation (\ref{MIS}), with initial iterate $u_0$ converges monotonically (non-decreasing) and uniformly towards to the solution ${u(x)}$ of (\ref{SBVP-1})--(\ref{SBVP-2}). Similarly $v_0$ as an initial iterate leads to a non-increasing sequences $\{v_n\}$ converging to a solution $v(x)$. Any solution $z(x)$ in $D$ satisfies
\begin{eqnarray*}
 u(x)\leq z(x) \leq  v(x).
\end{eqnarray*}
\end{theorem}
\textbf{Proof:} It is easy to show that (see \cite{AKV-NA-2011})
\begin{eqnarray*}
u_{0}\leq u_{1}\leq u_{2}\leq  \ldots \leq u_{n}\leq u_{n+1}\leq  \ldots \leq v_{n+1}\leq v_{n}\leq \ldots \leq v_{2}\leq v_{1}\leq v_{0}.
\end{eqnarray*}
So the sequences $u_{n}$ is monotonically non-decreasing and bounded above by $v_0$, similarly $v_{n}$ is non-increasing, respectively and bounded below by $u_0$. Hence by Dini's theorem they converges uniformly. Let ${u}(x)=\displaystyle\lim_{n\to\infty}u_{n}(x) $ and ${v}(x)=\displaystyle\lim_{n\to\infty}v_{n}(x)$.

Now by using Lemma \ref{Lemma2.3}, the solution $u_{n+1}$ of (\ref{MIS}) is given by
\begin{eqnarray*}
u_{n+1}= \frac{b ~ x^{\nu}J_{-\nu}\left(x \sqrt{\lambda }\right)}{J_{-\nu }\left(\sqrt{\lambda }\right)-\delta  \eta ^{\nu } J_{-\nu }\left(\eta  \sqrt{\lambda }\right)}- {\int_0}^1{G(x,t)t^{\alpha} (f(t,u_{n})-\lambda u_{n})dt}.
\end{eqnarray*}
Now by Lebesgue's dominated convergence theorem, we get
\begin{eqnarray*}
&&{u}(x)= \frac{b ~ x^{\nu}J_{-\nu}\left(x \sqrt{\lambda }\right)}{J_{-\nu }\left(\sqrt{\lambda }\right)-\delta  \eta ^{\nu } J_{-\nu }\left(\eta  \sqrt{\lambda }\right)}- {\int_0}^1{G(x,t) t^{\alpha} (f(t,{u})-\lambda {u})dt}.
\end{eqnarray*}
Which is the solution of nonlinear SBVPs (\ref{SBVP-1})--(\ref{SBVP-2}). Similar equation we can define for the sequence of lower solution also. It is easy to see that ${u}(x)\leq z(x)\leq {v}(x)$.

\subsubsection{Well ordered upper $(u_0)$ and lower $(v_0)$ solutions $(v_0\leq u_0).$}
Based on the sign of $\lambda$, we prove two Existence theorems; Theorem \ref{Theorem3.2} and Theorem \ref{Theorem3.3}. The proof of these theorems are similar to the proof of Theorem \ref{Theorem3.1}.
\begin{theorem}\label{Theorem3.2}
Assume that
\begin{itemize}
\item[$(W_1)$] there exist  $u_0$, $v_0$ in $C^2[0,1]$ such that $v_0\leq u_0$, where $u_0$ satisfies (\ref{US-1})--(\ref{US-2}) and $v_0$ satisfies (\ref{LS-1})--(\ref{LS-2});
\item[$(W_2)$] the function $f: D_1 \rightarrow R$ is continuous on $D_1 :=\{(x,y) \in [0,1] \times R :v_0 \leq y \leq u_0\}$;
\item[$(W_3)$] there exist $M\geq 0$ such that for all $(x,y),(x,w)\in D_1$,
\begin{eqnarray*}
y\leq w \Longrightarrow  f(x,w)- f(x,y) \geq M_2(w - y),
\end{eqnarray*}
\item[$(W_4)$] there exist a constant $\lambda>0$ such that $M_2-\lambda\geq0$ and $H_0$ or $H_2$ holds;
\end{itemize}
then the nonlinear three point SBVP (\ref{SBVP-1})--(\ref{SBVP-2}) has at least one solution in the region $D_1$. Sequence $\{u_n\}$ generated by equation (\ref{MIS}), with initial iterate $u_0$ converges monotonically (non-increasing) and uniformly towards a solution ${u(x)}$ of (\ref{SBVP-1})--(\ref{SBVP-2}). Similarly $v_0$ as an initial iterate leads to a non-decreasing sequences $\{v_n\}$ converging to a solution $v(x)$. Any solution $z(x)$ in $D_1$ satisfies
\begin{eqnarray*}
v(x)\leq  z(x) \leq u(x).
\end{eqnarray*}
\end{theorem}
\begin{theorem}\label{Theorem3.3}
Assume that
\begin{itemize}
\item[$(W'_1)$] there exist  $u_0$, $v_0$  in $C^2[0,1]$ such that $v_0\leq u_0$, where $u_0$ satisfies (\ref{US-1})--(\ref{US-2}) and $v_0$ satisfies (\ref{LS-1})--(\ref{LS-2});
\item[$(W'_2)$] the function $f: D_2 \rightarrow R$ is continuous on $D_2:=\{(x,y) \in [0,1] \times R :v_0 \leq y \leq u_0\}$;
\item[$(W'_3)$] there exist $ M_3 \geq 0$ such that for all $(x,\widetilde{y}),(x,\widetilde{w})\in D_2$
\begin{eqnarray*}
\widetilde{y}\leq \widetilde{w} \Longrightarrow  f(x,\widetilde{w})- f(x,\widetilde{y}) \geq  -M(\widetilde{w} - \widetilde{y})
\end{eqnarray*}
\item[$(W'_4)$] there exist a constant $\lambda<0$ such that $ M_3+\lambda \leq 0$ and $(H_0')$ holds;
\end{itemize}
then the nonlinear three point SBVP (\ref{SBVP-1})--(\ref{SBVP-2}) has at least one solution in the region $D_2$. Sequences $\{u_n\}$ generated by equation (\ref{MIS}), with initial iterate $u_0$ converges monotonically (non-increasing) and uniformly towards a solution ${u(x)}$ of (\ref{SBVP-1})--(\ref{SBVP-2}). Similarly $v_0$ as an initial iterate leads to a non-decreasing sequences $\{v_n\}$ converging to a solution $v(x)$. Any solution $ z(x)$ in $D_2$ satisfies
\begin{eqnarray*}
v(x)\leq z(x) \leq u(x).
\end{eqnarray*}
\end{theorem}

\subsection{Uniqueness of nonlinear three point SBVP}
\begin{theorem}
Let $f(x,y)$ be continuous on $D$ (or $D_1$ or $D_2$) and there exist a constant $M_{\lambda}\geq0$ such that
\begin{eqnarray*}
&& f(x,u)- f(x,v) \leq M_{\lambda}(u - v),
\end{eqnarray*}
and $M_{\lambda}<\lambda_1$, where $\lambda_1 \leq \min{\left\{j^{2}_{\nu,1},y^{2}_{\nu,1},i^{2}_{-\nu,1},k^{2}_{\nu,1}\right\}}.$ Then the nonlinear three point SBVP (\ref{SBVP-1})--(\ref{SBVP-2}) has a unique solution.
\end{theorem}
\textbf{Proof:} Suppose $u(x)$ and $v(x)$ be any two solutions of (\ref{SBVP-1})--(\ref{SBVP-2}) then we get
\begin{eqnarray*}
&&-(x^{\alpha} (u-v)')'= x^{\alpha}[f(x,u)-f(x,v)],\\
&&(u-v)'(0)=0,~~~~(u-v)(1)= \delta (u-v)(\eta),
\end{eqnarray*}
which gives
\begin{eqnarray*}
&&-(x^{\alpha} (u-v)')'- M_{\lambda} x^{\alpha}(u-v)\leq 0,\\
&& (u-v)'(0)=0,~~~~(u-v)(1)= \delta (u-v)(\eta).
\end{eqnarray*}
Since by the Maximum and Anti-maximum principles $(b=0)$, whenever $M_{\lambda}<\lambda_1$, we get $u-v\leq 0$ or $u-v\geq 0$ (i.e., $u\leq v$ or $u\geq v$) for different class of $\alpha$. Similarly by changing the role of $u$ and $v$, we get $u\geq v$ or $u\leq v$. Hence $u \equiv v$. Therefor the solution of the (\ref{SBVP-1})--(\ref{SBVP-2}) is unique.

\section{Numerical illustrations}\label{Example}
We present here some numerical examples to validate our existence results which is derived in the Theorem \ref{Theorem3.1}, Theorem \ref{Theorem3.2} and Theorem \ref{Theorem3.3}.
\subsection{Reverse ordered upper and lower solution}
The following examples validate the result of Theorem \ref{Theorem3.1}, and gives a range of $\lambda$ for that we can generate two monotone sequences which converge to the solution of nonlinear SBVP.
\begin{example}
Consider the nonlinear three point SBVP
\begin{eqnarray}
\label{Ex4.1.1}&& -y''(x)-\frac{\alpha}{x} y'(x)=\frac{\alpha~(e^{y}-1)-x}{4},\\
\label{Ex4.2.1}&&~~y'(0)=0,~~~~~~ y(1)=3 y\left(\frac{1}{7}\right).
\end{eqnarray}
\end{example}
Where $f(x,y)=\frac{\alpha~(e^{y}-1)-x}{4},~~\delta=3,~\eta =\frac{1}{7}$ and $\alpha$ satisfies $(H_1)$ or $(H_3)$. In this problem we choose lower and upper solutions as $v_0= 1$ and $u_0=-1$, where $v_0\geq u_0$, i.e., it is reverse ordered case. Here nonlinear term satisfies all assumptions for Theorem \ref{Theorem3.1} and Lipschitz constant is $M_1=\frac{e\alpha}{4}$. Now we can find out a range for $\lambda\in (\frac{e\alpha}{4},~~j_{\nu,1}^2)$ such that $(H_1)$ or $(H_3)$ are true.
\subsection{Well ordered upper and lower solution}
The following examples validate the results of Theorem \ref{Theorem3.2} and Theorem \ref{Theorem3.3}. On the basis of sign of $``\lambda"$, we divide this subsection into the following two parts.
\subsubsection{When $\lambda >0$}
\begin{example}
Consider the nonlinear three point SBVP
\begin{eqnarray}
\label{Ex4.1.2}&& -y''(x)-\frac{\alpha}{x} y'(x)=de^{y},\\
\label{Ex4.2.2}&&~~y'(0)=0,~~~~~~ y(1)=\frac{1}{2} y\left(\frac{1}{3}\right).
\end{eqnarray}
\end{example}
where $f(x,y)=de^{y},~~d=\frac{2(1+\alpha)e^{-\frac{2}{3}}}{3}~~\delta=\frac{1}{2},~\eta =\frac{1}{3}$ and $\alpha$ satisfies $(H_0)$ or $(H_2)$. In this problem we choose lower and upper solutions as $v_0= 0$ and $u_0=\frac{2-x^2}{3}$, where $v_0\leq u_0$, i.e., it is well ordered case. The nonlinear term satisfies all assumptions for Theorem \ref{Theorem3.2} and Lipschitz constant is $M_2=d$. Now we can find out a range for $\lambda>0$ such that the conditions $M_2-\lambda \geq 0$,~$(H_0)$ or $(H_2)$ are true.
\subsubsection{When $\lambda <0$}
\begin{example}
Consider the nonlinear three point SBVP
\begin{eqnarray}
\label{Ex4.1.3}&& -y''(x)-\frac{\alpha}{x} y'(x)=1-2y^3,\\
\label{Ex4.2.3}&&~~y'(0)=0,~~~~~~ y(1)=\frac{1}{3} y\left(\frac{1}{4}\right).
\end{eqnarray}
\end{example}
Here $f(x,y)=1-2y^3,~~\delta=\frac{1}{3},~\eta =\frac{1}{4}$ and $\alpha$ satisfies $(H'_0)$. In this problem we choose lower and upper solutions as $v_0= -1$ and $u_0=1$, where $v_0\leq u_0$, i.e., it is well ordered case. The nonlinear term satisfies all assumptions for Theorem \ref{Theorem3.3}, and Lipschitz constant is $M_3=6$. Now we can find out a range for $\lambda<0$ such that the conditions $M_3+\lambda \leq 0$ and $(H'_0)$ are true.
\begin{example}
Consider the nonlinear three point SBVP
\begin{eqnarray}
\label{Ex4.1.4}&& -y''(x)-\frac{3}{x} y'(x)=1-7y^2,\\
\label{Ex4.2.4}&&~~y'(0)=0,~~~~~~ y(1)=2.2 y\left(\frac{1}{5}\right).
\end{eqnarray}
\end{example}
Here $f(x,y)=1-7y^3,~~\delta=2.2,~\eta =\frac{1}{5}$ and $\alpha$ satisfies $(H'_0)$. In this problem we choose lower and upper solutions as $v_0= 0$ and $u_0=\frac{7}{4}+\frac{5}{2}x^2$ i.e., it is well ordered case. The nonlinear term satisfies all assumptions for Theorem \ref{Theorem3.3}, and Lipschitz constant is $M_3=\frac{119}{2}$. Now we can find out a range for $\lambda<0$ such that the conditions $M_3+\lambda \leq 0$ and $(H'_0)$ are true.

\end{document}